\documentclass[10pt,a4paper]{article}
\usepackage[latin1]{inputenc}
\usepackage{anysize}
\usepackage{graphicx}
\marginsize{2.5cm}{2.5cm}{2cm}{1cm}
\usepackage{amssymb,amsmath,amscd,amsfonts,amsthm,mathrsfs}
\newcommand{\bb}[1]{\mathbb{#1}}

\numberwithin{equation}{section}
\newtheorem{theorem}{Theorem}
\newtheorem{proposition}[theorem]{Proposition}
\newtheorem{lemma}[theorem]{Lemma}
\newtheorem{example}[theorem]{Example}
\newtheorem{corollary}[theorem]{Corollary}
\usepackage{url}
\begin{document}
\title{A Note on Combinatorial Derivation}
\author{Joshua Erde \thanks{DPMMS, University of Cambridge}}
\maketitle
\begin{abstract}
Given an infinite group $G$ and a subset $A$ of $G$ we let $\Delta(A) = \{g \in G \,:\, |gA \cap A| =\infty\}$ (this is sometimes called the \emph{combinatorial derivation} of $A$). A subset $A$ of $G$ is called \emph{large} if there exists a finite subset $F$ of $G$ such that $FA=G$. We show that given a large set $X$, and a decomposition $X=A_1 \cup \ldots \cup A_n$, there must exist an $i$ such that $\Delta(A_i)$ is large. This answers a question of Protasov. We also answer a number of related questions of Protasov. 
\end{abstract}
\section{Introduction}
For a subset $A$ of an infinite group $G$ we denote

$$\Delta(A) = \{g \in G \,:\, |gA \cap A| = \infty\}.$$

This is sometimes called the \emph{combinatorial derivation} of $A$. We note that $\Delta(A)$ is a subset of $AA^{-1}$, the difference set of $A$. It can sometimes be useful to consider $\Delta(A)$ as the elements that appear in $AA^{-1}$ `with infinite multiplicity'. In \cite{P2011} Protasov analysed a series of results on the subset combinatorics of groups (see the survey \cite{S2011}) with reference to the function $\Delta$, and asked a number of questions. In this note we present answers to some of those questions. \\
\ \\
A subset $A$ of $G$ is said to be \cite{LP2009}:
\begin{itemize}
\item \emph{large} if there exists a finite subset $F$ of $G$ such that $FA=G$;
\item \emph{$\Delta$-large} if there exists a finite subset $F$ of $G$ such that $F \Delta (A) = G$.
\end{itemize}
Protasov asked \cite{P2011}:\\
\ \\
{\bf Question A.}
\emph{Is every large subset of an arbitrary infinite group $G$ $\Delta$-large?}\\
\ \\
{\bf Question B.}
\emph{Does there exist a function $f:\bb{N} \rightarrow \bb{N}$ such that, for any group $G$ and any partition $G=A_1 \cup \ldots \cup A_n$, there exists an $i$ and a subset $F$ of $G$ such that $G=F \Delta(A_i)$ and $|F| \leq f(n)$?}\\
\ \\
{\bf Question C.}
\emph{Let $G$ be an infinite group. Given $G=A_1 \cup \ldots \cup A_n$, such that $A_i=A_i^{-1}$ and $e \in A_i$ for $i \in \{1,\ldots,n\}$, does there exist an $i$ and an infinite subset $X$ of $G$ such that $X \subseteq A_i$ and $\Delta(X) \subseteq A_i$?}\\
\ \\
A subset $A$ of $G$ is said to be \emph{sparse} if for every infinite subset $X$ of $G$, there exists a non-empty finite subset $F$ of $X$ such that $\bigcap_{g \in F} gA$ is finite. A subset $A$ of $G$ is said to be \emph{$\nabla$-thin} if either $A$ is finite, or there exists an $n \in \bb{N}$ such that $\Delta^n(A) = \{e\}$, where $\Delta^n$ denotes the iterated application of $\Delta$. Protasov also asked:\\
\ \\
{\bf Question D.}
\emph{Is every $\nabla$-thin subset of a group $G$ sparse?}\\
\ \\
In Section $2$ we present answers to the four questions posed by Protasov. We say a subset $A$ of $G$ is \emph{cofinite} if there exists a finite subset $H$ of $G$ such that $A = G \setminus H$. Our main result is:

\begin{theorem}\label{t:union}
Let $G$ be an infinite group. Given a subset $X$ of $G$ such that there exists a finite subset $F$ of $G$ such that $FX$ is cofinite, and a decomposition $X=A_1 \cup \ldots \cup A_n$, then there exists an $i$ and a subset $F'$ of $G$ such that $|F'| \leq |F|(|F|+1)^{2^{n-1}-1}$ and $F'\Delta(A_i) = G$.
\end{theorem}

This provides a positive answer to Questions A and B. We also show:

\begin{theorem}\label{t:delta}
Let $G$ be an infinite group. Given a subset $A$ of $G$ and a countable subset $X$ of $\Delta(A)$ such that $X=X^{-1}$ and $e \in X$, there exists a subset $Y$ of $A$ such that $\Delta(Y) = X$.
\end{theorem}

Since for all infinite sets $A$ we have that $e \in \Delta(A)$, this provides a positive answer to Question C. Finally we also show:

\begin{theorem}\label{t:sparse}
Let $G$ be an infinite group, $A$ a subset of $G$. Then if $A$ is $\nabla$-thin, then $A$ is sparse.
\end{theorem} 
\section{Results}
We will start by considering Question A. For a subset $A$ of $G$ and a finite subset $F$, it is easy to see that $\Delta(FA) = F\Delta(A)F^{-1}$. Therefore for abelian groups it is apparent that if $A$ is large, with say $FA=G$, we have that $F^{-1}F\Delta(A) = \Delta(G) = G$, and so $A$ is $\Delta$-large. 
\begin{lemma}\label{l:large}
Let $G$ be an infinite group, $A$ a subset of $G$. If there exists a finite subset $F$ of $G$ such that $FA$ is cofinite, then $F \Delta(A) = G$.
\end{lemma}
\proof
We would like to find some finite set $X = \{x_1, \ldots , x_k\}$, such that the set of translates $\{x_1A , x_2A \ldots , x_kA\}$ has the property that, for any $g \in G$, we must have $|gA \cap x_iA| = \infty$ for some $i$. Then, for all $g \in G$ we would have that $|x_i^{-1}gA \cap A| = \infty$ for some $i$ and so $x_i^{-1}g \in \Delta(A)$ and so $g \in X\Delta(A)$. Therefore we could conclude that $X\Delta(A) = G$. Let $F=\{f_1,\ldots,f_k\}$.\\
\ \\
Since $FA$ is cofinite, there exists some finite subset $H$ of $G$ such that $f_1A \cup \ldots \cup f_kA = FA = G \setminus H$. Therefore we see that for any $g \in G$ there must exist an $i$ such that $|gA \cap f_iA| = \infty$. Therefore $F$ satisfies the property above, and so $F\Delta(A) = G$.
\qed \\
\ \\
We note that Question A follows from Lemma \ref{l:large} as a simple corollary.
\begin{corollary}\label{c:large}
Let $G$ be an infinite group, $A$ a subset of $G$. Then if $A$ is large, then $A$ is $\Delta$-large. \qed \end{corollary}

There are various concepts of ``small" for subsets of groups. A subset $A$ of $G$ is said to be \cite{LP2009}:
\begin{itemize}
\item \emph{small} if $(G \setminus A) \cap L$ is large for every large subset $L$ of $G$;
\item \emph{P-small} if there exists an injective sequence $(g_n)_1^{\infty}$ in $G$ such that $g_i A \cap g_j A = \phi$ for all $i,j$; 
\item \emph{almost P-small} if there exists an injective sequence $(g_n)_1^{\infty}$ in $G$ such that $|g_i A \cap g_j A| < \infty$ for all $i,j$;
\item \emph{weakly P-small} if for every $n \in \bb{N}$ there exists distinct elements $g_1,g_2, \ldots , g_n$ in $G$ such that $g_i A \cap g_j A = \phi$ for all $i,j$.
\end{itemize}
Protasov also asked:\\
 \ \\
{\bf Question E.}
\emph{Is every nonsmall subset of an arbitrary infinite group $G$ $\Delta$-large?}\\
\ \\ 
We can also use the method in Lemma \ref{l:large} to show that sets which are not almost P-small are $\Delta$-large.
\begin{theorem}\label{t:P-small}
Let $G$ be an infinite group, $A$ a subset of $G$. Then if $A$ is not almost P-small, then $A$ is $\Delta$-large.
\end{theorem}
\proof
Take a maximal set $F=\{f_1, \ldots , f_k\}$ such that $|f_i A \cap f_j A| < \infty$ for all $i,j$. Such a set exists and is finite since $A$ is not almost P-small. Then, for all $g \in G$ we must have that $|gA \cap f_i A| = \infty$ for some $i$, since $F$ is maximal. Hence $f_i^{-1}g \in \Delta(A)$ and so $G = F \Delta (A)$.
\qed \\
\ \\
We note however that there do exist sets $A$ which are not weakly P-small (and so also not P-small), but which are still not $\Delta$-large. 
\begin{example}\label{e:P-small}
Consider the group $(\bb{Z},+)$. Let $A=\{10^n\,,\,10^n+n \,:\, n \in \bb{N} \}$. Clearly any translate of $A$ has non-empty intersection with $A$, and so $A$ cannot be weakly P-small. However $\Delta(A) = \{0\}$ since each difference only appears a finite number of times in $A$.
\end{example} 
It remains to show that sets which are not small are $\Delta$-large. We will be able to show this using ideas that are used in answering Question B.\\
\ \\
To motivate the proof we first consider the case $n=2$. Given a decomposition of $G$ into two sets $A \cup B$ what does it mean if $\Delta(A) \neq G$? Well in that case we have some $g \in G$, $g \not\in \Delta(A)$. Therefore there are only a finite number of $h \in G$ such that the group elements $h$ and $g^{-1}h$ are both members of $A$, since each such $h$ is in $gA \cap A$. Therefore there is some finite subset $H$ of $G$ such that for all $h \in G \setminus H$, either $h \in B$ or $g^{-1}h \in B$. But then we have that $\{ e , g \}B = B \cup g B = G \setminus H$ and so, by Lemma \ref{l:large}, $\{ e  ,g \}\Delta(B) = G$. This idea motivates the following lemma which will be key to answering Question B.

\begin{lemma}\label{l:union}
Let $G$ be an infinite group. Let $X$ be a subset of $G$ such that there exists finite subsets $F,H_1$ of $G$ such that $FX = G \setminus H_1$. Then given a decomposition of $X$ into two sets $X=A \cup B$, either $F\Delta(A) = G$ or there exists $g \in G$ and a finite subset $H_2$ of $X$  such that $(gF \cup \{e\}) B = X \setminus H_2$.
\end{lemma}
\proof
 Let $F = \{f_1, \ldots , f_k\}$. If $F\Delta(A) \neq G$, then there exists $g \in G$, $g \not\in F\Delta(A)$, that is, $f_i^{-1} g  \not\in \Delta(A)$ for $i=1,\ldots,k$. Now, as before, for each $i$ there are only finitely many $h \in X$ such that both $h$ and $f_i^{-1} g h  \in A$. Also we claim that there are only finitely many $h$ such that none of the group elements $f_1 ^{-1} g h, \ldots , f_k ^{-1} g h $ are in $X$. Indeed, since if $F^{-1} g h \cap X = \phi$ then we have that $gh \cap FX = gh \cap G \setminus H_1 = \phi$ and so $h \in g^{-1} H_1$.\\
\ \\
Therefore we have that there exists some finite subset $H_2$ of $X$ such that for all $h \in X \setminus H_2$ and for all $i$, no pair $h$, $f_i^{-1} g h$ are both in $A$, and at least one of the group elements $f_i^{-1} g h$ is in $X$. Therefore we have that $B \cup g^{-1}FB = X \setminus H_2$.
\qed \\
\ \\
We note at this point that this lemma allows us to settle the final part of Question E, whether or not a subset $A$ of $G$ which is not small, must be $\Delta$-large.
\begin{corollary}\label{c:small}
Let $G$ be an infinite group, $A$ a subset of $G$. Then if $A$ is not small, then $A$ is $\Delta$-large.
\end{corollary}
\proof
If $A$ is not small then there exists a large set $L$ such that $(G \setminus A) \cap L$ is not large. Without loss of generality let us assume that $A \subset L$. Then $L = (L \setminus A) \cup A$, and there exists a finite subset $F$ of $G$ such that $FL = G$. Therefore, by Lemma \ref{l:union}, if $F \Delta(A) \neq G$, then there exists $g \in G$ and a finite subset $H_2$ of $L$ such that $(gF \cup \{e\}) (L \setminus A) = L \setminus H_2$. Therefore there is some finite subset $H_3$ of $G$ such that $F(gF \cup \{e\}) (L \setminus A) = G \setminus H_3$. However it is then clear that there exists some finite subset $F'$ of $G$ such that $F'(L \setminus A) = G$, however by assumption $L \setminus A$ was not large, and so $F \Delta(A) = G$. Therefore $A$ is $\Delta$-large.
\qed \\
\ \\
We can apply Lemma \ref{l:union} inductively to prove the main result in the note. \\
\ \\
\emph{Proof of Theorem \ref{t:union}.}
We induct on $n$. The case $n=1$ follows from Lemma \ref{l:large}. Given that the result holds for all $k < n$, and given a subset $X$ of $G$ and finite subsets $F,H_1$ of $G$ such that $FX = G \setminus H_1$ and a decomposition $X=A_1 \cup \ldots \cup A_n$. We have, by Lemma \ref{l:union}, that either $F\Delta(A_1) = G$, or there exists $g \in G$ and a finite subset $H_2$ of $G$ such that $(gF \cup \{e\}) (A_2 \cup \ldots \cup A_n) = X \setminus H_2$. We apply the induction hypothesis to the set $Y=(A_2 \cup \ldots \cup A_n)$. Note that there exists a finite subset $H_3$ of $G$ such that $F(gF \cup \{e\})Y = G \setminus H_3$. Thus we have that there exists an $i$ and a subset $F'$ of $G$ such that

\begin{align*}
|F'| &\leq |F|(|F|+1)\big(|F|(|F|+1) + 1\big)^{2^{n-2}- 1} \leq |F|(|F|+1)\big((|F|+1)^2\big)^{2^{n-2}-1}\\
 &\leq |F|(|F|+1)(|F|+1)^{2.2^{n-2}-2} \leq |F|(|F|+1)^{2^{n-1}-1}
\end{align*}

and $F' \Delta(A_i) = G$.
\qed \\
\ \\
It is known \cite{BP2003} that given a decomposition $G = A_1 \cup \ldots \cup A_n$, there exists an $i$ and a subset $F$ of $G$ such that $|F| \leq 2^{2^{n-1}-1}$ and $F(A_iA_i^{-1})= G$. As $\Delta(A)$ is a subset of $AA^{-1}$ the following Corollary strengthens this result.
\begin{corollary}\label{c:union}
Let $G$ be an infinite group. Given a decomposition $G=A_1 \cup \ldots \cup A_n$ then there exists an $i$ and a subset $F$ of $G$ such that $|F| \leq 2^{2^{n-1}-1}$ and $F\Delta(A_i) = G$.
\end{corollary}
\proof
We apply Theorem \ref{t:union} with $X=G$ and $F= \{e\}$ to get the result stated. \qed\\
\ \\

Theorem \ref{t:union} says that if we decompose any large set into a finite number of pieces, at least one of the parts must be $\Delta$-large. However there do exist $\Delta$-large sets which decompose into two sets which are not $\Delta$-large.

\begin{example}
Consider the group $(\bb{Z},+)$. Let $A = \{10^n : n \in \bb{N}\}$ and let $B = \{10^1 + 1\} \cup \{10^2 + 1\} \cup \{10^3 + 2\} \cup \{10^4 + 1\} \cup \{10^5 + 2\} \cup \{10^6 + 3\} \cup \{10^7 + 1\} \ldots$. Then if $X = A \cup B$ we see immediately that $\Delta(X) = \bb{Z}$, but $\Delta(A) = \Delta(B)=\phi$.
\end{example}

There also exist decompositions of large sets into a finite number of sets, none of which are large.

\begin{example}
Consider the free group on $2$ elements, $F(a,b)$. If we denote by $aSb$ the set of reduced words in $F(a,b)$ that start with $a$ and end with $b$, then it is clear that $aSb$ is not large. Indeed no finite set of translates can contain the words $a^n$ for all $n$. However we can decompose $F(a,b)$ as

$$F(a,b) \setminus \{e\} = \bigcup_{x,y \in \{a,a^{-1},b,b^{-1}\}} xSy,$$
none of which are large.
\end{example}
We now consider Question C. In \cite{P2011} it is shown that for all infinite groups $G$, and all subsets $A$ of $G$ such that $A=A^{-1}$ and $e \in A$, there exists some subset $X$ of $G$ such that $\Delta(X) = A$. Using a similar construction we are able to prove Theorem \ref{t:delta}.\\
\ \\
\emph{Proof of Theorem \ref{t:delta}.}
Let $X=\{x_1,x_2, \ldots\}$ and let $Z_0 = Y_0 = \phi$. We define a new sequence. For all $n \in \bb{N}$ let
\begin{equation*}
w_i = x_{i-\frac{n(n-1)}{2}} \,\,\,\,\text{ for }\,\,\,\,\, \frac{n(n-1)}{2} + 1 \leq i \leq \frac{n(n+1)}{2}.
\end{equation*}
That is, $w_1 = x_1$, $w_2=x_1,w_3 = x_2$, $w_4 = x_1, w_5 = x_2, w_6 = x_3$, $w_7 = x_1 \ldots$. Our plan is to add pairs of elements to $Y$ such that each pair has difference $w_i$, but introduces no new differences with the elements already in $Y$. Given $Y_{i-1} = \bigcup_{j=0} ^{i-1} Z_j$ we want to inductively find $Z_{i} = \{ z_i, w_i z_i\} \subset A$ such that $(Z_i Y_{i-1}^{-1} \cup Y_{i-1}Z_i^{-1}) \cap Y_{i-1} Y_{i-1}^{-1} = \phi$. Equivalently $Z_i$ needs to avoid the finite set $Y_{i-1}Y_{i-1}^{-1}Y_{i-1}$.  This is always possible since $w_i \in \Delta(A)$ and so the number of such pairs is infinite. We let $Y = \bigcup_{i=1}^{\infty} Y_i$ and see that $\Delta(Y)=X$ as claimed.
\qed \\
\ \\
We note at this point that, perhaps surprisingly, it is necessary for the subset $X$ in Theorem \ref{t:delta} to be countable.
\begin{proposition}
There exists a group $G$, a subset $A$ of $G$, and a subset $Y$ of $\Delta(A)$ such that there does not exist any subset $X$ of $A$ with $\Delta(X) = Y$.
\end{proposition}
\proof
Let $\kappa = \big( 2^{\aleph_0}\big)^+$ (so assuming CH we would have that $\kappa = \aleph_2$), and let $\alpha$ be the initial ordinal of cardinality $\kappa$. Consider the group $G = (\bb{Z}_2)^{\alpha}$. That is, $G$ is the direct product of $\kappa$ copies of $\bb{Z}_2$. Let $x_i$, for $i \leq \alpha$, be be the element which is $1$ in the $i$th copy of $\bb{Z}_2$ and $0$ elsewhere. Consider the set

$$ A = \{x_n + x_i \,:\, 1\leq n < \omega \,,\, \omega \leq i \leq \alpha\} \cup \{x_n \,:\, 1\leq n < \omega\}.$$

We see that the set $X = \{x_i\,:\, \omega \leq i \leq \alpha\}$ is a subset of $\Delta(A)$, however we claim that there does not exist any subset $Y$ of $A$ such that $\Delta(Y) = X$. Indeed, suppose such a subset exists. Let $i$ be such that $\omega \leq i \leq \alpha$. Since $x_i \in \Delta(Y)$ we have that the set of $n < \omega$ such that both $x_n$ and $x_n + x_i$ are in $Y$ is infinite. Let us call this set $L_i$. Since $2^{\aleph_0} < \kappa$ we must have that there exist $i,j \leq \alpha$ such that $L_i = L_j$. But then we also have that $x_i + x_j \in \Delta(Y)$, contradicting our initial assumption. \qed \\

This phenomenon arises since in the definition of the function $\Delta$ we only require that the intersection $|gA \cap A|$ is infinite. If instead we were to consider a more general function

$$\Delta_{|G|} (A) = \{g \in G \,:\, |gA \cap A| = |G|\},$$ 

an analogous version of Theorem \ref{t:delta} could be proved, by the same argument, for sets of larger cardinality. Indeed much of the work in this note, including more general versions of Lemma \ref{l:large} and Theorem \ref{t:union}, can be easily adapted to state results in this framework. Given Theorem \ref{t:delta} we can answer Question C, albeit in a slightly trivial way.

\begin{corollary}\label{c:thin}
Let $G$ be an infinite group, $G=A_1 \cup \ldots \cup A_n$, $A_i=A_i^{-1}$, $e \in A_i$ for $i \in \{1,\ldots,n\}$. Then there exists an $i$ and an infinite subset $X$ of $G$ such that $X \subseteq A_i$ and $\Delta(X) \subseteq A_i$.
\end{corollary}

\proof
By Theorem \ref{t:delta} we have that, within any infinite set $A$ there exists a subset $X$ of $A$ with $\Delta(X) = \{e\}$. At least one of the $A_i$ must be infinite, and therefore such an $X$ satisfies the statement.
\qed

\ \\
Finally we turn to Question D.\\
\ \\
\emph{Proof of Theorem \ref{t:sparse}.}
We will prove that every set which is not sparse is not $\nabla$-thin. Given $A \subset G$ which is not sparse then there exists some infinite subset $X$ of $G$ such that for every finite subset $F$ of $X$ $\bigcap_{g \in F} gA$ is finite. In particular for every pair $g_i,g_j \in X$ we have that $|g_iA \cap g_j A| = \infty$ and so $X^{-1}X \subset \Delta(A)$.\\
\ \\
However for an infinite set $X$ we claim that $X^{-1}X \subset \Delta(X^{-1}X)$. Indeed given $g_i^{-1}g_j \in X^{-1}X$ we have that 

$$|g_i^{-1}g_j X^{-1}X \cap X^{-1}X| = |g_j X^{-1}X \cap g_i X^{-1}X| \geq |X \cap X| = |X| = \infty,$$

and so $g_i^{-1}g_j \in \Delta(X^{-1}X)$. Therefore if $A$ is not sparse we have that $X^{-1}X \subset \Delta^n(A)$ for all $n \in \bb{N}$, and so $\Delta^n(A) \neq \phi$ for any $n$.
\qed \\
\ \\
However, there exist $2$-sparse sets which are not $\nabla$-thin. 

\begin{example}
Consider the group $(\bb{Z},+)$. Let

\begin{align*}
A = \{0\} &\cup \{1\} \cup \{10\} \cup \{11\} \cup \{10^2\} \cup \{10^2 + 3\} \cup \{10^3\} \cup \{10^3 +1\} \\
 &\cup \{10^4\} \cup \{10^4 + 3 \} \cup \{10^5\} \cup \{10^5 + 5 \} \cup \{10^6\} \cup \{10^6 + 1 \} \cup \ldots 
\end{align*}
Then we see that $\Delta(A) = \{0\} \cup \{2n + 1 \,:\, n \in \bb{Z}\}$, and so $\Delta^n(A) = \{2n \,:\,n \in \bb{Z}\}$ for all $n \geq 2$. However, given any infinite subset $X$ of $G$ we must have two numbers $a$ and $b$ in $X$ whose difference is even. Then $aA \cap bA$ is finite, since there are only a finite number of even differences in $A$ less than any given number.
\end{example}
Finally we mention an open problem. Corollary \ref{c:union} says that if an infinite group $G$ is split into a finite number of sets, then one of those sets must be $\Delta$-large. For groups of larger cardinality can we prove similar results? For example:\\
\ \\
{\bf Question F.}
\emph{Let $G$ be an infinite group, with $|G|= \kappa$. Given $\mu < \kappa$, $|I|=\mu$, and a decomposition $G = \bigcup_{i \in I} A_i$, does there exist an $i \in I$ and a finite subset $F$ of $G$ such that $F\Delta(A_i)=G$?} 
\bibliography{A-note-on-Combinatorial-Derivation}
\bibliographystyle{plain}
\end{document}